\documentclass[12pt,a4paper]{article}
\usepackage[utf8]{inputenc}
\usepackage[T1]{fontenc}
\usepackage[english]{babel}
\usepackage{amsmath, amssymb, amsfonts}
\usepackage{booktabs}
\usepackage{geometry}
\usepackage{xcolor}
\usepackage{float}
\usepackage{tikz}
\usepackage{times}
\usepackage{url}
\usepackage{hyperref}
\hypersetup{
	colorlinks=true,
	linkcolor=blue,
	urlcolor=blue,
	citecolor=blue
}

\title{\textbf{An Old Babylonian coefficient, its origin and impact on our understanding of measures on circles, including the radian measure.}}
\author{Jens Kleb \\
	\small ORCID: 0000-0002-3304-2253}
\date{\small {1.November 2025}}

\begin{document}
	
	\maketitle
	
	\begin{abstract}
		This study reconstructs the origin of a constant, here called $\Xi$ (Xi), as a primary scaling factor in Old Babylonian mathematics and astronomy. $\Xi$ arises from the practical necessity of precise measurements on the sky or a circle, through the harmonization of length-measure systems. The analysis of the Nippur measure (with its famous cubit) and the Gudea measure shows that $\Xi = 375/360$ represents the ratio of these established Old Babylonian measure systems. As a precision factor for circumference calculations, it remained in use until today. In Ptolemy's work, we find a slightly refined value of $\Xi = 377/360$. A further refinement of this coefficient led to our modern $\pi$, which still incorporates the two Old Babylonian components of a demonstrably two-stage calculation and refinement process. The accuracy increased by only 0.5\% compared to the first ratio. This factor, attested on several Old Babylonian cuneiform tablets including those from Susa, demonstrates the profound understanding of sexagesimal logic. The relative sexagesimal notation (60 = 1 = 1/60) enabled the universal application of $\Xi$ and its reciprocal for highly accurate calculations of arc-length on circular segments. This investigation leads ultimately to a surprising consequence: the modern radian measure is a direct descendant of this Old Babylonian coefficient.
	\end{abstract}
	
	\section{Introduction}
	
	The history of the mathematical constant $\pi$ is closely linked to the development of ancient mathematics and astronomy. Fundamental research by Otto Neugebauer, as well as the comprehensive account by C. J. Scriba and P. Schreiber, has significantly shaped our understanding of Babylonian mathematics \cite{Neugebauer1945,Neugebauer1954,ScribaSchreiber2005}. R.C.A. Rottländer's metrological studies over decades support, that our current understanding of $\pi$ as the primary circular constant may be a logical consequence of earlier considerations — and that another quantity, here denoted by $\Xi$ (Xi), reflects these earlier considerations historically and conceptually \cite{Muroi, Rottlaender1979,Kleb2023a,Kleb2023b}.
	
	The highly precise astronomical and geometrical tasks of the Babylonians required the exact determination of circular segments. Scriba and Schreiber describe Babylonian mathematics as ``essentially an algebraically and number-theoretically oriented science'' that operated with concrete numbers \cite[p. 16]{ScribaSchreiber2005}. The central question was: how does one obtain the curved arc length $b$ of a sector from the straight line (segment width or chord length)?
	
	\subsection{Commemorative Notes on Two Pioneers of Cuneiform Mathematics}
	
	Before presenting the analysis, a brief tribute to two scholars whose work on mathematical cuneiform tablets laid essential foundations is appropriate.
	
	\subsubsection*{Otto Neugebauer (1899--1990)}
	With meticulous philological rigour and mathematical clarity, Neugebauer deciphered Babylonian mathematics as an autonomous, highly logical discipline. His edition and interpretation of mathematical cuneiform tablets liberated ancient algebra and geometry from the shadow of Greek tradition—fundamentally reshaping the history of exact sciences.
	
	\subsubsection*{Evert Marie Bruins (1909--1990)}
	Bruins reconstructed the procedural core of Mesopotamian tablet mathematics through detailed analyses of calculation techniques. His studies, revealed sophisticated number-theoretic practices among scribes, making fraction arithmetic and area computations in cuneiform contexts systematically accessible for the first time.
	
	\subsection{The Sexagesimal Foundation: 60 = 1}
	
	The basis of $\Xi$ is the Babylonian sexagesimal system, in which numbers are conceived as place values to base 60. Neugebauer demonstrated the paramount importance of this number system for Babylonian arithmetic and astronomy \cite[pp. 20--29]{Neugebauer1969}. 
	
	In this sexagesimal thinking, the hexagon served as the basis of the circle. Initially, the radius and 1/6 of the circumference were considered equal. Through the here introduced $\xi$ (Xi) factor, the circumference portion was scaled (blowed up) to a near-real arc-length.
	
	A fundamental peculiarity is the equivalence of 60 and 1—depending on the place, the number ``1" represents either one unit, 60 units, or their higher powers, as well as $1/60$. This place-value logic transfers to geometry: a circle of radius \textit{60} (length units) is, in sexagesimal thinking, identical to a circle of radius \textit{1} \cite{Neugebauer1945,Neugebauer1954,Kleb2022}.
	
	\textbf{Important remark:} Sexagesimal numbers are written here in \textit{italic} with places separated by ``;", conforming to original cuneiform notation.
	
	\section{The Metrological Origin of Xi}
	
	The crucial insight lies in Old Babylonian metrology. The \textbf{Nippur measure system} and the \textbf{Gudea measure system} stood in a ratio of 360/375 = \textbf{0.96} (sexagesimal \textit{57;36}) \cite{Rottlaender1979, Kleb2023b}. Its reciprocal number, \textit{1;2;30}, reveals the emergence of Xi:
	
	To ensure that exactly \textbf{360 units of the Nippur system} remained on the precisely calculated circumference, the diameter unit had to be reduced by a factor of 0.96. The practical implementation was the invention of an adapted measure system during Gudea's reign: 120 Gudea units in diameter yielded 375 Gudea units on the corresponding circumference—while maintaining the 360 Nippur length units on this circumference as a undoubted basis for partitioning the circle.
	The relationship between these two systems of length, was not only discussed by works of Rottländer and the author (\cite{Rottlaenderwebsite,Kleb2023a}) , but can also be measured by anyone from the proportions of the famous ``Nippur-cubit" in the museum of Istanbul and from the precisely recorded length engraved on Gudea's statue, in the Louvre. Artifacts from Gudea´s reign and Lifetime in addition give us the reason why he officially has established this reform of measures too.
	
	\begin{center}
		\textbf{375 / 360 = $\Xi$ = \textit{1;2;30} and its sexagesimal reciprocal \textit{57;36}}.
	\end{center}
	
	Thus:
	\[
	\Xi_1 = \frac{375}{360} = \frac{\textit{6;15}}{\textit{6}} = \textit{1;2;30} \approx 1.041666\ldots
	\]
	It is important to know that the factor $\Xi_1 = 375/360$ can be reduced, resulting in $25/24$ and its reciprocal, $24/25$. Although mathematically clearly equivalent, this reducing obscures the historical origin: The numerator $375$ derives from the number of Gudea units on the circumference at a diameter of \textit{2} = $120$ ($4 \times 30= 4$ cubits), and the denominator $360$ represents the traditional division of the circle into 360 parts, inherited from the Nippur system of measurement. Retaining this ratio $375/360$ in the following paragraphs therefore preserves the metrological context and makes the transformation between the two coexisting systems of length visible. This also makes clearer and more understandable that the same, only refined, factor is contained in the subsequent refinement in Ptolemy's time $377/360$ and also in our present-day $\pi = 3 \times 376.9911.../360$.\\
	
	This value, attested for the time of Gudea (ca. 2100 BCE), represents the oldest known form of $\Xi$ \textit{1;2;30} (therefore here called $\Xi_1$) \cite{Kleb2023b}.
	
	\subsection{Definition and Meaning of $\Xi$}
	
	The symbol and name, within this work established as $\Xi$, (from Greek $\alpha$$\upsilon$$\Xi$$\eta$$\sigma$$\iota$$\varsigma$ -``auXesis," meaning "growth" or "increase") denotes the coefficient for stretching from original 360 parts (hexagon-parts) to 375 diameter units (as real arc-length) on the circumference. 
	It was determined by the Babylonians dividing, the number of diameter units (the Gudea unit shrinked from Nippur) required for the now precisely calculated circumference of the circle, by 360 ``old" units/parts on circumference a circle partitioning that was very essential and ``unchangeable" in Babylonian thinking as it is even until today. Like the Diameter, of \textit{``2"} $= 120$ Units equal to 4 Cubits, also the Partitioning of circles, coming from 12 $\times$ 30 or 6 $\times$ 60 into 360 parts was mentioned and used, several times on cuneiform texts in sexagesimal notation as simply \textit{``6"}.
	
	To illustrate how the same sexagesimal number \textit{1;2;30} represents different absolute values depending on context (the Babylonian floating-point interpretation), the following equivalences hold:
	
	\vspace{0.2cm}
	\noindent
	\begin{tabular}{@{}p{0.23\textwidth} p{0.33\textwidth} p{0.39\textwidth}@{}}
		\toprule
		\textbf{Sexagesimal Value} & \textbf{Interpretation (Babylonian)} & \textbf{Modern, just  refined Equivalent} \\
		\midrule
		\textit{1;2;30} & Stretch factor for the whole circle and its parts, in decimal especially for 6; 60 or 360 (\textit{6} hexagon parts $\rightarrow$ \textit{6;15} and so on) & $\pi / 3 \approx 1.047197\ldots$ \\ 
		\textit{1;2;30} & Stretch factor for 1/360 of the circle (\textit{1} part $\rightarrow$ \textit{1;2;30} parts) & $\pi / 180 \approx 0.0174533\ldots$ \\
		\textit{57;36} & its Reciprocal: shrink factor from arc to hexagon-part & $180 / \pi \approx 57.29578\ldots$ \\
		\bottomrule
	\end{tabular}
	\vspace{0.2cm}
	\begin{center}
		{\LARGE $\Xi$} = {\Large \textit{1;2;30}} , reciprocal: {\Large \textit{57;36}}
	\end{center}
	In sexagesimal place-value system, $\Xi$ is dimensionless and can apply to 1, 60, $1/60$, or 3600 depending only on context \cite{Kleb2021,Kleb2023b}.\\
	\\
	\textbf{Historical and mathematical evolution:}
	\begin{itemize}
		\item \textbf{First step:} \textit{\large{3}} — the simple Babylonian ratio between diameter and circumference simplified by a Hexagon and its 6 isosceles triangles\cite{CDLI2}.
		\item \textbf{Second step:} ${\large \Xi}$ — the coefficient for precise calculations on circle and the circumference (late 3rd to early 2nd millennium BC).
		\[
		{\large \Xi_1 = \frac{\textit{1;2;30}}{\textit{1}} = \frac{\pi_{susa}}{\textit{3}} = \frac{\textit{3;7,30}}{\textit{3}}}
		\]
		\item \textbf{Third step:} ${\large \pi}$ — the product of both factors before, if precision was necessary.
	\end{itemize}
	\begin{center}
		\textit{3} $\cdot$ $\Xi_1$ = \textit{3;7;30} = $\pi_{susa} \:\:(3.125$)
	\end{center}
	
	\subsection{Selected Cuneiform Sources}
		Consistent nomenclature is established by cross-referencing Louvre (Sb) and Yale (YBC) collections, the very worth-full ORACC Database by Eleanor Robson \& Laith M. Hussein \cite{ORACC}, CDLI identifiers checked at Oct. 2025:
	\begin{table}[H]
		\centering
		\begin{tabular}{@{}p{3.0cm}p{3.6cm}p{6.2cm}p{2.0cm}@{}}
			\toprule
			\textbf{Inventory No.} & \textbf{Publication} & \textbf{Physical Content} & \textbf{CDLI Ref.} \\ \midrule
			\textbf{Sb 13089}\cite{CDLI3} & TMS 3 (MDP 34, 3) & Coefficients (reciprocal of $\Xi_1$ \textit{57;36}) & P254835 \\
			\textbf{Sb 13088}\cite{CDLI2} & TMS 2 (MDP 34, 2) & regular Hexagon Drawing & P254827 \\
			\textbf{Sb 13090}\cite{CDLI1} & TMS 1 (MDP 34, 1) & Mathematical sketch ($\Xi_1$ implicit) & P254816 \\
			\textbf{YBC 7243}\cite{CDLI4} & ORACC \& MCT & Coefficient of a log \textit{4;48} & P255045 \\
			\textbf{YBC 5022}\cite{CDLI5} & Neugebauer/Sachs MCT & Coefficients (reciprocal of $\Xi_1$ \textit{57;36}) & P255026 \\ 
			\textbf{YBC 8600}\cite{Neugebauer1945} & Neugebauer/Sachs MCT & Mathematical calculation of a log by \textit{4;48}, reciproc \textit{12;30} & P255065 \\ \bottomrule
		\end{tabular}
	\end{table} 
	
	\section{Confirmation on Cuneiform Tablets}
	
	\subsection{The Value \textit{57;36} on Cuneiform Tablets}
	
	The term UB (Akkadian tubqum), traditionally translated as "sector" or "corner," serves in Babylonian mathematics not as an abstract value for angular measurements, but as a designation for a base segment, or count of vertices of a polygon.
	
	The coefficient \textit{57;36} (reciprocal of \textit{1;2;30}) functions as a geometric transformation factor: it calculates the ratio between a linear portion of a polygon and the corresponding arc length on a circle. This value "stretches" the simple hexagon to a more precise circular geometry.
	
	For a regular polygon with N vertices:
	\begin{center}
		\textbf{(\textit{6/N}) $\cdot$ \textit{1;2;30} $\cdot$ r = arc length}
	\end{center}
	\begin{center}
		\textbf{{\Large $\frac{ \textit{6} \cdot \textit{r}}{\textit{arc length} \cdot \textit{57;36}}$ = N = $\frac{\textit{6;15} \cdot \textit{r}}{\textit{arc length}}$}}
	\end{center}
	
	This coefficient implicitly leads to $\pi_{susa} = 3.125$. Neugebauer showed that a hexagon with side length 1 (60) has perimeter \textit{6} (360), while the circumscribing circle has perimeter \textit{6;15} (375): \textit{6 / 6;15 = 57;36} \cite{Neugebauer1969}. Hofmann (1964), Beckmann (1970), and others also calculated a $\pi_{susa}$ of 3.125 using this tablet and line in various ways \cite{BIB-Orient, Beckmann}.
	
	A second important mention appears on YBC 5022, lines 65-66, with the term \textit{``ša eb-lu-ú-um"} ("of the rope"). The rope symbolizes both the stretched connection between vertices (as in a hexagon) and the curved arc like on a loop. The coefficient \textit{57;36} represents the ratio of the entire circumference to its hexagon, as well as its individual segments.
	\begin{center}
		\textbf{\textit{Arc length} $\cdot$ \textit{0;57;36} = hexagon portion} \:\:\: or also \:\:\:\: \textbf{\textit{hexagon portion} $\cdot$ \textit{1;2;30} = arc length} (this can be multiplied now by the given radius to complete the calculation.)
		for instance: 1/25th of circle at radius 1 is a hexagon portion of \textit{14;24 $\cdot$ 1;2;30 $\cdot$ r(1)  = 15}\\ With our modern further refined coefficient and in decimal numbers and also at radius 1 , the arc-length is 0.251327... .
	\end{center} 
	The formulas are valid for any fractal or segment belonging to the same radius and corner or sector width, the so called ``UB" of the first mentioned tablet (equal to an arbitrary part of this circle or polygon).
	
	\subsection{Terminological Note: UB and rope}
	
	The interpretation of the coefficients \textit{57;36} and \textit{1;2;30} as circumference-related stretching factors is not a modern imposition but is grounded in the original edition and translation of the cuneiform tablets.
	
	For TMS 3 (Sb 13089), the first editors, Bruins in a preliminary article and Bruins \& Rutten in detail, described the coefficient in line 30 explicitly as ``la constante du cycle (cercle plus parfait)'' — the constant of the (more perfect) circle \cite{BruinsRutten1961,Bruins1950}. This phrasing leaves no ambiguity: the coefficient pertains to the circumference of an idealized circle beyond the simple hexagon approximation.
	
	The term ``UB'' for the same tablet and line is a newer translation by E.Robson, who meticulously re-analyzed and listed the cuneiform texts in ORACC. This has now ``corner'',``edge" or ``sector'' \cite{ORACC}. While this might seem at first glance very different, but it isn´t, in Babylonian mathematical practice it denotes a segment of the perimeter corresponding to a distance between vertices or a sector of circles. This is consistent with previous interpretations and also clarifies the connections through Robson's adapted more precise translation.
	
	For YBC 5022, line 65–66, the term ``rope'' (Akkadian \textit{eblu}) appears. Robson again translates this in a geometrically neutral way, but the context — a list of coefficients placed under the authority of the goddess Nisaba, who demands precise measures — leads to a very plausible interpretation, that the ``rope'', symbolizes twice at same time: The curved arc (the rope laid around the circle like a loop) , and on the other side, the straight segments of a hexagon perimeter (the rope stretched between vertices). Therefore, the coefficient \textit{57;36} scripted direct before the mention of the goddess, underlines its high importance.  The ratio of the hexagon perimeter to the circumference of the circle \cite{Neugebauer1945, Neugebauer1969}.
	
	Thus, even the most cautious and critical modern translations do not dispute the essential point: the coefficients \textit{1;2;30} and \textit{57;36} are tied to the circle's circumference and its division into parts. An alternative interpretation as a pure area coefficient can not supported by the primary publications or by newer contextual analysis of the tablets.\\ Finally, it should be mentioned that this symbol was surely also chosen because the traditional unit of length for a rope: a Ninda was 12 cubits in length. A hexagon with radius of 2 has exactly this circumference (12). However, its surrounding circle was in real half a cubit longer. This results in: $12 / \textit{12;30} = \textit{57;36}$. This is the factor by which the rope (Ninda) that fits the perimeter of the hexagon, is shorter than the circumference of its circle, a really perfect synonym in Babylonian thinking. 
	
	\subsection{Intermediate Results}
	
	Both tablets, along with refined coefficients for circle area (i.e. \textit{5} to \textit{4;48}), demonstrate Babylonian refinement. Neugebauer and Sachs suspected such refinements as early as 1945 \cite[p. 59]{Neugebauer1945}. The factor \textit{4;48} appears as \textit{5 $\cdot$ 57;36 = 4;48}.
	
	The Babylonians divided circles into 6 $\times$ 60 parts. A full circumference had up to 21,600 parts (sexagesimal \textit{6}), with 375 to 22,500 diameter units (sexagesimal \textit{6;15}). The ratio 0.96/1 corresponds to \textit{57;36}, with reciprocal \textit{1;2;30} = $\Xi_1$.
	
	The mathematical cuneiform tablets from Susa contain coefficients but not directly a $\pi$-value of \textit{3;7;30} \cite{Neugebauer1945,BruinsRutten1961}. This $\pi$ follows necessarily from ``3" and $\Xi_1$ as a second refinement stage.
	
	Thus the evolution is clear:
	\begin{itemize}
		\item \textbf{Not:} They had $\pi_{susa} = 3.125$ and derived $\Xi_1 = \pi_{susa}/3$.
		\item \textbf{Rather:} They had \textit{3} and the refinement factor $\Xi_1$ \textit{1;2;30} with reciprocal \textit{57;36}, created from metrological necessity, and combined them when precision was needed.
	\end{itemize}
	\begin{center}
		\Large{$3 \cdot \Xi_1 = \pi_{susa} = 3.125$} \cite{BruinsRutten1961,Neugebauer1969,Kleb2023a,Kleb2023b}
	\end{center}
	The sophisticated use of modular formulas and  if necessary a refinement factor—such as the \textbf{$\Xi_1$ $\textit{1;02;30}$} and the reciprocal as ``shrink factor" for the high precision level \textbf{$\textit{57;36}$}—demonstrates a geometric use on terrestrial but of prime importance for astronomical calculations. \\
	\\
	In Eratosthenes' time and before, instead of the direct angle, the chord was determined as leg of a triangle, and this portion was then considered as part of a full circle divided into 360 or by this as an regular polygon with N-corners/edges.\\
	The Greek word ``$\mu$$o$$\iota$$\rho$$\alpha$$\iota$" is mentioned several times in ancient literature, like also the Roman ``Parcae" and``Pars", or the parts named in arabic scientific history ``Dakaica"(\cite{Neugebauer1962},p17) all corresponding to a division by 360, just as degrees correspond to one 360th of a circle and its smaller parts of 60th each in sexagesimal notation. Eratosthenes himself used the sunbeam, shadow length, and shadow-stick as a triangle to determine the angle of sun in Alexandria (relative to the zenith in Syene)\cite{PtolemyAlmagest, Kleb2023b}. His triangle's side ratio, according to the tradition of plausible measure instrument dimensions described here, which have a length of 4 cubits (4$\times$30 = 120), was probably $\approx$15.16/120 this leads to N=50 or 1/50th of a circle.
	\subsection{The Two-Stage Refinement Process}
	
	\begin{table}[H]
		\centering
		\small
		\caption{Formulary: Two-stage refinement process}
		\label{tab:Formulary}
		\begin{tabular}{@{}p{3.6cm}p{1.2cm}p{2.1cm}p{6.2cm}p{1.9cm}@{}}
			\toprule
			\textbf{Target Value} & & \textbf{Primary Source} & \textbf{Operational Path} & \textbf{Adapted Coefficient} \\ \midrule
			\textbf{Diameter} ($d$) & & & $d = (c_{Susa} \times \textit{\textbf{20}}) \times \textit{\textbf{0;57;36}}$ & $\textit{\textbf{19;12}}$\\
			\textbf{Circumference} ($c$) & & & $c_{Susa} = (d \times \textit{\textbf{3}}) \times \textit{\textbf{1;02;30}}$ & $\textit{\textbf{3;07;30}}$ \\
			\textbf{Vertices/Sector} ($UB$) & & TMS 3 & N = $\frac{\textit{\textbf{6}} \times d/2}{Arc\text{-}length \times \textit{\textbf{57;36}}}$ & $\textit{\textbf{6;15}}$ \\
			\textbf{Arc-Length} ($UB$) & & TMS 3 & $b = (\textit{\textbf{6}}/N \times d/2) \times \textit{\textbf{1;02;30}}$ & $\textit{\textbf{6;15}}$ \\
			\textbf{Hexagon-part} ($Rope$) & & YBC 5022 & $Hexagon\text{-}part = b \times \textit{\textbf{57;36}}$ & $\textit{\textbf{57;36}}$ \\
			\textbf{Arc-Length} ($Rope$) & & YBC 5022 & $b = Hexagon\text{-}part \times \textit{\textbf{1;02;30}}$ & $\textit{\textbf{1;02;30}}$ \\
			\textbf{Area} ($log$) & & YBC 7243 & $A_{Susa} = c_{Susa}^2 \times (\textit{\textbf{5}} \times \textit{\textbf{57;36}})$ & $\textit{\textbf{4;48}}$ \\
			\textbf{Diameter/Perimeter} ($log$)\cite{CDLI4} & & YBC 8600 & $3d = \textit{57;36} \times \sqrt{A_{Susa} \times \textit{\textbf{12}} \times \textit{\textbf{1;02;30}}}$ & $\textit{\textbf{12;30}}$ \\
			\midrule
			\textbf{Precision circle} & & & $\pi_{susa} = Hex. perimeter \textit{\textbf{``3"}} \times \: \textit{\textbf{1;02;30}}$ & $\textit{\textbf{1;02;30}}$ \\
			\midrule
			Stretch Factor & $\textit{1;02;30}$ && $A_{Susa} = A \times \textit{1;02;30}$\\
			Shrink Factor & $\textit{0;57;36}$ && $V = V_{Susa} \times \textit{57;36}$ \\
			\bottomrule
		\end{tabular}
	\end{table}
	
	\subsection{Interpreting the Evidence: Plausibility and Continuity}
	
	The cuneiform tablets cited in this study date primarily to the Old Babylonian period (ca. 1900–1600 BCE), while Gudea's metrological reform belongs to the late third millennium (ca. 2100 BCE). This temporal gap of several centuries raises the question whether the coefficient \textit{1;2;30}, reciproc \textit{57;36} were indeed inherited from Gudea's time or represent a later independent development. While absolute certainty is unattainable after four millennia, several converging arguments support the continuity hypothesis.
	
	First, the metrological relationship between the Nippur and Gudea systems ($360/375$) is archaeologically attested through preserved standards and is not dependent on the later tablets. The tablets thus employ a ratio that demonstrably existed centuries earlier.

	Second, the context of the coefficients on the tablets themselves points to a deliberate connection with the tradition of precise measures. On YBC 5022, line 66, the coefficient \textit{57;36} is directly followed by a reference to the goddess Nisaba, the patron of writing, accounts, and — as the Gudea cylinders make explicit — of divinely mandated measurement. Neugebauer interpreted the presence of such a reference in a list of coefficients as an indication that this text belonged to the genre of ``wisdom literature'' disseminated in Nisaba's name \cite{Neugebauer1945}. In this cultural context, mathematical precision was not merely practical but carried cosmic and religious significance, a theme already prominent in Gudea's inscriptions.
	
	Third, the very structure of our modern $\pi$ — defined as $3 \cdot (\pi/3)$ — still reflects the two-stage Babylonian approach. The refinement from $\Xi_1 = 1.041666\ldots$ to $\Xi_2 = 1.04722\ldots$ to the modern $\Xi = 1.047197\ldots$ increases precision by only 0.5\% over 4000 years, suggesting that the conceptual framework was established early and underwent only incremental numerical improvement. Had the Babylonians of the Old Babylonian period derived their coefficients independently, it would be a remarkable coincidence that they exactly preserved the Gudea ratio.
	
	Ultimately, historical reconstruction of ancient mathematics deals in probabilities, not proofs. What can be asserted is that a coherent tradition of circle measurement based on the ratio $375/360$ is traceable from the late third millennium through the Old Babylonian period to Ptolemy and beyond. The tablets from Susa and the Yale collection provide tangible evidence that this tradition was actively cultivated by scribes who understood the sexagesimal flexibility of their numbers and who, like Gudea before them, placed this knowledge under the authority of Nisaba.
	
	\subsection{On the Interpretation of Coefficients: Circumference vs. Area, and the Indirect Evidence for Precision}
	
	A recurring point of discussion in the secondary literature concerns whether the coefficients \textit{57;36} and \textit{1;2;30} should be understood as relating to the circumference or to the area of a circle. Bruins \& Rutten initially described the value \textit{57;36} in TMS 3, line 30 as ``la constante du cycle (cercle plus parfait)'' \cite{BruinsRutten1961}, explicitly linking it to the circle's perimeter. Neugebauer independently arrived at the same conclusion, demonstrating that for a hexagon with side length 1 (equal to 60), the perimeter of the circumscribed circle is \textit{6;15} (375), yielding the ratio \textit{6 / 6;15 = 57;36} \cite{Neugebauer1969}. 
	
	It is true that other tablets, such as YBC 8600 and YBC 5022, contain area coefficients like \textit{4;48} (which is $\textit{5} \times \textit{57;36}$) and \textit{12;30} (the reciprocal of $\textit{4;48}$). These area coefficients are derived from the circumference coefficients by multiplication with the constant 5, which in Babylonian mathematics served as the conversion factor from circumference to area for a circle when using the approximation $\pi = 3$ (i.e., $c^2/12$). The introduction of the refinement factor $\Xi_1$ accordingly transforms the area coefficient from $\textit{5}$ to $\textit{5} \times \textit{57;36} = \textit{4;48}$, as noted by Neugebauer \cite[pp. 57-59]{Neugebauer1945}. 
	
	Thus, the primary coefficients \textit{1;2;30} and \textit{57;36} are fundamentally related to the circumference. The area coefficients that appear on the same tablets are secondary, derived by multiplication with the constant 5 (or its reciprocal). This distinction is crucial: the stretching factor $\Xi$ operates on linear measures (the circumference and its arc segments) and only indirectly affects area calculations through the established Babylonian formula $A = c^2 \times \textit{5}$ (or its refined version $A = c^2 \times \textit{4;48}$). 
	
	\paragraph{The sagitta problems like on BM 85194.}
	Further support comes from geometrical problems involving the \textit{sagitta} (arrow, i.e., the height of a circular segment) and inscribed rectangles. Both Höyrup \cite{Hoyrup2002} and Friberg \cite{Friberg2005}, who have been perfectly analyzing ancient Babylonian problems for decades, demonstrate in some problems the relationship between the chord, the sagitta and the diameter in Babylonian mathematical thinking. These calculations using the formula $s = r - \sqrt{r^2 - (c/2)^2}$ (also see \cite{Kleb2022,Kleb2021}) . However, this formula is geometrically exact only for a true circle. Otherwise the Babylonians clear used the hexagon as a simplification like real circle by the knowing, that the inscribed rectangle will only fit if there is a curved arc around. The consequence is: Such calculations implicitly require a correct relationship between the linear dimensions of the circle and its curved segments. The solution methods often employ coefficients that implicitly incorporate a more precise circle constant than the simple hexagon approximation. Moreover, these problems are not limited to the hexagon; they work with arbitrary chords and segments, which presupposes a general theory of circle measurement consistent with the two-stage refinement process described here.  
	
	\paragraph{Indirect evidence from astronomy and geometry.} 
	The necessity of a precise circumference coefficient beyond the simple hexagon is in addition independently confirmed by the astronomical records transmitted us from Babylonia. The exact orbital periods of the Moon, planets, and the calculation of synodic and sidereal cycles require a circle divided into more than the six segments of a hexagon; they demand a division into at least 360 parts and their fractals by a proportional factor that accurately relates linear measures to arc lengths. The Venus tablets of Ammisaduqa, the lunar ephemerides, and the later development of the System A and B lunar theories all presuppose a circle with sufficient precision to accumulate small segments without accumulating unacceptable error \cite{Neugebauer1969}. Such precision would be impossible if the circumference were calculated merely as $3 \times d$, the simple hexagon approximation, especially if the observations could be done only in small segments of the sky, which have to be combined exactly to a circle.
	
	\paragraph{The refinement process in practical application.}
	Even though the lower accuracy level (level 1) was sufficient for most calculations and materials, it turns out that the more precise level 2 was sometimes used for expensive and important materials such as fine-grained substances requiring precise volume calculations, or simply because of their weight. This is evident not only in the case of clearly round logs (log, YBC 5022; YBC 8600), whose value and weight were important, but also in the case of precise coefficients for grain piles (YBC 5022), where their volume, i.e., the actual quantity, was paramount. This also certainly played a role in determining the actual load-bearing capacity of a boat (IM 49949) or a jug too.  This is just a short selection of such materials and geometric forms, that could easily be refined/precised if needed by using the optional two-stage calculation process. 
	The emergence of such refined coefficients corresponded precisely to the goal of the Gudea reform. The existence of such area and volume coefficients derived from the perimeter factor does not call into question the interpretation of \textit{57;36} as a perimeter-related scaling factor; rather, it confirms the modular, two-stage refinement process described in this study, in which a single coefficient can be applied to different geometric contexts.
	
	\section{Restoring the Inner Order: The Dream of Gudea and the Mandate of Nisaba}
	
	The metrological reform under Gudea of Lagash (ca. 2144–2124 BCE) was not merely a practical adjustment but a deliberate act of cosmic restoration. This is vividly documented in the famous Gudea cylinders (Cylinders A and B), inscribed in Sumerian on two large clay cylinders discovered at Girsu (modern Tello) and now housed in the Louvre Museum \cite{CDLI_Gudea, Edzard1997, Roemer2010}. These texts, among the longest known Sumerian literary compositions, recount in over 1,400 lines how Gudea received a divine command to rebuild the E-ninnu temple for the city god Ningirsu \cite{Rey2020}.
	
	The narrative begins with a crisis: the annual floods of the Tigris had failed, famine threatened the land, and cosmic order was disrupted \cite{Bewer1919}. In a nocturnal vision, Gudea saw his master, Lord Ninĝirsu, who spoke to him of his house and its building \cite{ETCSL_c2.1.7}. Yet the message remained enigmatic. Seeking clarity, Gudea travelled by canal to the temple of the goddess Nanše, the divine dream-interpreter. Along the way, he made offerings at shrines, preparing himself for revelation \cite{ETCSL_c2.1.7}.
	
	The dream itself contained potent symbols. A giant figure—with wings, a crown, and two lions—commanded him to build. Two other figures appeared: \textbf{a woman holding a golden stylus, studying a clay tablet on which the starry heaven was depicted}, and a hero carrying a lapis lazuli tablet with the plan of a house \cite{WorldHistory_Nisaba}. When Nanše interpreted the dream, she explained that the giant was her brother Ninĝirsu; the woman with the golden stylus was \textbf{Nisaba}, the goddess of writing, accounts, and—crucially—\textbf{the keeper of divine measurements}. According to the interpretation, Nisaba was directing Gudea to lay out the temple \textbf{astronomically aligned with the “holy stars”} \cite{Verderame2020, LiuJian2016}. The hero was Nindub, the architect-god, surveying the plan \cite{WorldHistory_Nisaba}.
	
	Nisaba’s role in this vision is foundational. As the Sumerian goddess of writing and scribes, she was also the guardian of accounts, measures, and architectural plans \cite{WorldHistory_Nisaba}. Her symbol was the stylus; her domain included mathematics, astronomy, and the precise recording of earthly and celestial dimensions. In the early Dynastic period, cylinder seals depict her in association with construction, particularly of monuments and temples \cite{WorldHistory_Nisaba}. The goddess’s demand for precise measurement was not arbitrary: it was the prerequisite for restoring what the Sumerians called \textbf{\textit{ni-si-sa}}—the “inner order” or “cosmic justice” that Gudea was divinely mandated to re-establish \cite{Rottlaender1979, Kleb2023b}.
	
	Gudea’s response was to “restore the inner order of the world” through a dual act: a metrological reform that harmonized two measure systems (the Nippur and Gudea cubits), and the construction of a temple whose dimensions reflected celestial patterns. The ratio $375/360 = \textit{1;2;30}$ thus becomes more than a mathematical curiosity; it is a material expression of a cosmic principle, aligning terrestrial measures with divine will.
	
	The result of this restoration, as described in Cylinder B, was a temporary return to an ideal state—a “golden age” that accompanied the completion of the temple. The text states:
	
	\begin{quote}
		“On the day that the King entered the temple, for seven days the maidservant was equal to her mistress, the slave and the master walked side by side; in his city the mighty and the lowly lay side by side; on the wicked tongue the (bad) words were changed (into good ones) … the sun let shine forth righteousness, Babbar trampled unrighteousness underfoot” (Cylinder B 17:18–18:11) \cite{Bewer1919, Kleb2023b}.
	\end{quote}
	
	This vision of social harmony—where oppression ceases, justice prevails, and abundance flows—was understood as the direct consequence of restoring proper measure. The goddess Nisaba, who in her earlier form was a grain deity, had become the guarantor of that order through the precision of writing, counting, and building \cite{WorldHistory_Nisaba}.
	
	The differential usage of the stretching factor $\Xi_1 = \textit{1;2;30}$ and its reciprocal $\textit{57;36}$ in cuneiform tablets (such as TMS 3 and YBC 5022) demonstrates that this metrological precision was not a one-time event but a scalable system. It allowed Gudea—and later Babylonian scribes—to toggle between the practical geometry of the hexagon (stage 1) and the more precise curved reality of the circle (stage 2) whenever accuracy was demanded by ritual, astronomical, or administrative tasks \cite{Kleb2023a, Kleb2023b}.
	
	The sophisticated use of modular formulas and the refinement factor thus reflects a worldview in which mathematics, metrology, and divine order were inseparable. What we might see as an abstract coefficient was, for the Babylonians, a tool for cosmic restoration—a means to make the earthly temple mirror the heavenly plan, under the watchful guidance of Nisaba, the goddess who “opened the mouth for seven reeds” and held the lapis lazuli tablet of the stars \cite{WorldHistory_Nisaba}.
	
	\begin{figure}[htbp]
		\centering
		\begin{tikzpicture}[scale=1.4]
			\draw[thick] (0,0) circle (2cm);
			
			\foreach \i in {0,1,2,3,4,5} {
				\coordinate (P\i) at (60*\i:2cm);
			}
			
			\draw[thick, blue] (P0) -- (P1) -- (P2) -- (P3) -- (P4) -- (P5) -- cycle;
			
			\foreach \i in {0,1,2,3,4,5} {
				\draw[red, thick] (0,0) -- (P\i);
			}
			
			\draw[<->, thick] (0, -2.2cm) -- (0, 2.2cm);
			\draw[<->, thick] (-0.6cm, 1cm) -- ( 0cm, 1.2cm);
			\draw[<->, thick] (0.0cm, -1.3cm) -- ( 0.6cm, -1.1cm);
			
			\node at (-0.1,2.35) [above right] {\textcolor{red}{Diameter $d = 2r=2 \cdot 60$}};
			\node at (0.2, 1.6) [below left] {30 parts};
			\node at (1.1, -1.3) [below left] {30 parts};
			\node at (P0) [right] {Vertex = N};
			\node at (-0.5,2.30) [above left] {\textcolor{blue}{Hexagon side = $r$}};
			\node at (1.2,-1.5) [below right] {Circumference $C_{circle} = \textit{3} \cdot \textit{1;2;30} \cdot r$};
			
			\node at (3.2,1.2) [above right, align=left] {%
				Sexagesimal basis: $2r = \textit{2}\:\:  (120)$ parts \\
				$C_{hex} \:\:\:= \color{blue} \textit{6}$\:\:\:\:\:\: \color{blue}(360 parts) \\
				$C_{circle} = \color{red}\textit{6;15}$ \color{red}(375 parts) \\
				$\Xi = \textit{1;2;30}$ (stretch factor)%
			};
		\end{tikzpicture}
		\caption{Regular hexagon inscribed in a circle. The radius is divided into 60 parts (sexagesimal basis). The hexagon side equals the radius. The refinement factor $\Xi_1 = \textit{1;2;30}$ stretches the hexagon perimeter of 6 (360 parts) to the precise circular circumference of $\textit{6;15}$ (375 parts).}
		\label{fig:hexagon_circle}
	\end{figure}
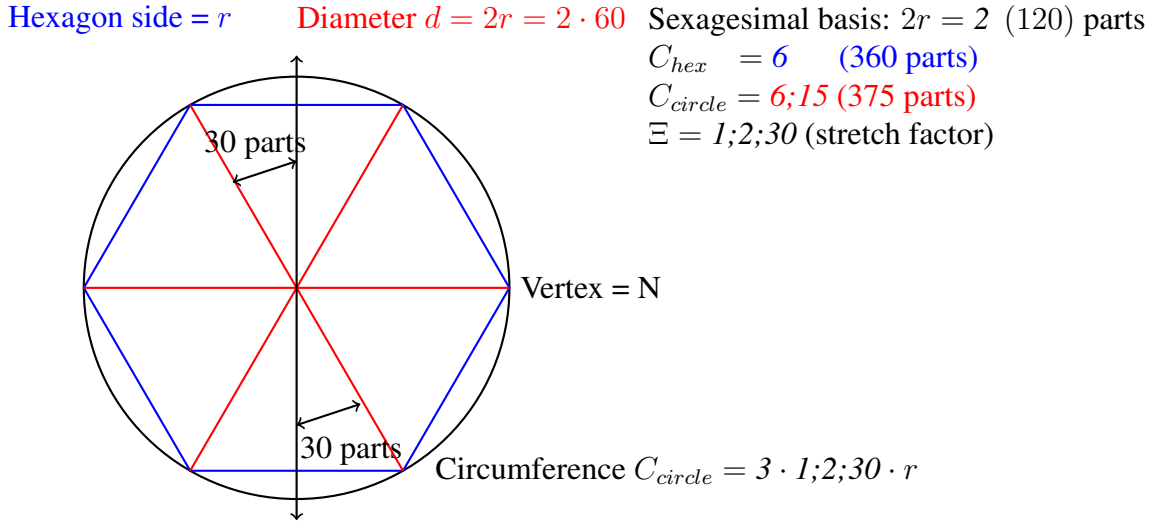
	
	\section{Refinement at Ptolemy's Time}
	
	The refinement from $\Xi_1$ ($\textit{1;2;30}=375/360$) to $\Xi_2$ ($\textit{1;2;50} = 377/360 = 1.04722...$) resulted from more precise calculations. Ptolemy's \emph{Almagest} consistently uses $\Xi_2$ as a practical factor \cite{PtolemyAlmagest} within his tables of chords.\\
	The period between Eratosthenes and Ptolemy, however, was marked by numerous changes, and some erroneous data clearly reflect these important shifts in calculation and data documentation \cite{Kleb2023c}. As mentioned in Historia I Swiat tome 13 (pp.230,\cite{Kleb2023b}), the data collected by Ptolemy and published in the \emph{Almagest} seems to contain both old and new datasets, sometimes mixed and misunderstood therefore. It seems inconceivable that Ptolemy would have provided calculations, that differed significantly from the results of his contemporaries. Ptolemy's data on the Moon's apparent diameter, when corrected from a chord to angle measure, otherwise align very good with given Chaldean averages (0° 31' 51"):
	
	\begin{table}[h]
		\centering
		\begin{tabular}{|l|c|c|c|}
			\hline
			& \textbf{min} & \textbf{max} & \textbf{averaged} \\
			\hline
			Original (Ptolemy) & \textit{31;20} & \textit{35;20} & \textit{33;20} \\
			After conversion to degrees & 0°29'55" & 0°33'45" & 0°31'50" \\
			\hline
		\end{tabular}
		\caption{Ptolemy's lunar data corrected from chords to degrees}
		\label{tab:modern}
	\end{table}
	It is in addition important to note that, the use of \textit{1;2;50} (377/360) marks a significant turning point, as this number lacked a direct reciprocal in the sexagesimal numbering system.
	For completeness: Ptolemy's $\pi$ value was by this refinement: \textit{3 $\times$ 1;2;50} = \textit{3;8;30} (= 377/120 = 3.141$\bar{6}$).
	
	\section{Comparison with Modern Radian Measure}
	
	\subsection{Deconstructing the Anachronism: Sexagesimal Floating Point}
	
	Right from the start, it must be stated that the following analysis is not an anachronistic equation of merely similar-looking mathematical quantities. Critics may claim that identifying $\textit{1;2;30}$ or its reciprocal $\textit{57;36}$ as an Old Babylonian kind of ``radian'' is anachronistic, because the Babylonians lacked the concept of $1^\circ$ or $rad$ as we understand them today. However, this objection fails to account for several crucial facts, just three of them will mentioned here now:
	
	First, it is beyond dispute that the Babylonians divided circle into at least 360 parts — coming from a segmentation of the sky or zodiac and also time — a tradition that around 2000 years later was adopted in our modern degree system. 
	
	Second, the sexagesimal place-value system is inherently \textbf{position-agnostic}: the same numerical notation $\textit{1;2;30}$ can represent $1.04166\ldots$ (the stretching factor for the whole circle) or $0.017361\ldots$ (the stretching factor for a single part and its fractals), depending solely on context. This is a typically Babylonian way of dealing with numbers, which makes it possible to perform several operations with one and the same sequence of numbers. 
	
	Third, the cuneiform tablets mentioned on this article, implicit or explicitly associate the coefficients $\textit{1;2;30}$ or its reciprocal $\textit{57;36}$ with circular geometry and figures.\\
	\\
	The sexagesimal system, which we still use today for time and angles, was therefore not an obstacle, but a very useful and efficient tool. That these facts and indications of its application, even in addition to the documented ancient Babylonian scientific and astronomical discoveries, are insufficient to prove its use in the manner described here is highly improbable, if not completely impossible, given the Babylonians' complex understanding of mathematics and geometry. This mathematical overview and efficiency should not lead us to doubt, but rather to have the greatest respect for our ancestors.
	
	\subsection{Direct Comparison: $\Xi$ and the Modern Radian}
	
	The modern radian measure defines the radian as the portion of a circle where the arc length equals the radius. Converting degrees to radians uses the factor $\pi/180$:
	
	\begin{center}
		Modern: $b = r \cdot (\pi/180) \cdot \alpha_{deg}$
	\end{center}
	
	The Susa tablet formula for the same arc length is:
	
	\begin{center}
		Babylonian: $b = r \cdot \Xi \cdot \text{hexagon-parts}$
	\end{center}
	
	Thus, for any angle measured in degrees ($\alpha_{deg}$) and in hexagon-parts (which are functionally identical to degrees), we have:
	
	\begin{center}
		$\alpha_{rad} = \alpha_{deg} \cdot \Xi$ \quad and therefore \quad $\Xi / 60 = \pi / 180$.
	\end{center}
	Furthermore, this is the decimal representation used today; in sexagesimal notation, division by 60 is not necessary, so the same coefficient is by this view also $\Xi = \pi /180$.
	Both coefficients are conceptually consistent and differ only in their degree of refinement. In other words, both represent the same underlying mathematical-geometric concept, with the modern version being a further refinement of the original Babylonian coefficient. 
	
	\subsection{The Evidence from Ptolemy's Chord Tables}
	
	The continuity of this conceptual framework is confirmed by Ptolemy's work around 150 CE. His \emph{Almagest} contains chord tables that implicitly use the refined value $\Xi_2 = \textit{1;2;50}$. This is exactly the same type of coefficient as the Old Babylonian $\Xi_1$, only slightly more precise. Ptolemy did not reinvent it; he refined the Babylonian stretching factor within the same sexagesimal framework. But, especially we can see Ptolemy´s $\Xi_2$ at the first line in its smallest value as a 60th of just one part of circle. This is $377/21,600 = 377/360 \cdot 60$. The fact that $\Xi_2$ lacks a perfect reciprocal (unlike $\Xi_1$) suggests that the original Old Babylonian design was deliberately optimized for arithmetic elegance, while later refinements prioritized geometric accuracy.
	
	\subsection{Numerical Comparison}
	
	On base 60, the number $\textit{1;02;30}$ is simultaneously:
	
	\begin{enumerate}
		\item The linear stretching factor from hexagon to circle perimeter: $1 + 2/60 + 30/3600 = 1.041\bar{6}$.
		\item The 1-part-to-arc factor (i.e., the radian for $1^\circ$): $1/60 + 2/3600 + 30/216000 = 0.017361\ldots$.
	\end{enumerate}
	
	This value ($0.017361\ldots$) is, aside from the modern refinement ($\pi/180 \approx 0.017453\ldots$), conceptually fully comparable in origin and evolutionary development to the radian for $1^\circ$. The only difference is that the Babylonian system was superior in its \textbf{arithmetic elegance} (finite, regular numbers with exact reciprocals), while our $\pi$ is more precise but irrational.
	
	\subsection{A modern or an Old Babylonian invention?}
	
	What we call the ``radian measure'' today is by this knowing not really a modern invention. It is the direct descendant of the Old Babylonian stretching factor $\Xi$, seen at this point in the context of 1 part of 360, (sexagesimal no difference, but in our decimal thinking one sixtieth of the initial decimal $\Xi$) just refined over four millennia from $\Xi_1 = 0.017361\ldots$, through $\Xi_2 = 0.0174537\ldots$ to the modern $\Xi = 0.0174533\ldots$. The Babylonians lacked the abstract concept of an angle measured in radians, but they possessed the functional equivalent: a proportional coefficient that related the number of parts (degrees) to the arc length. The sexagesimal system's place-value agnosticism allowed the same number to serve as both the radian for $60^\circ$ and the radian for $1^\circ$ — a flexibility that our decimal system cannot replicate without additional notation. The claim of anachronism thus collapses when one recognizes that mathematical concepts can exist in functional form long before they receive their modern names and definitions.

	\section{Conclusion}
	
	The analysis of the Nippur and Gudea measure systems has shown that the coefficient:\newline $\Xi = \textit{1;2;30} \:\: (375/360$) originated as a metrological necessity: it harmonized two coexisting length systems, which were well attested by several artifacts, while preserving the traditional division of a circle into 360 parts by the Babylonians. This coefficient, attested in cuneiform tablets from Susa and elsewhere (TMS 3, YBC 5022, YBC 8600 ...), served as an universal stretching factor that transformed the hexagon as a synonym for a simple circle (stage 1, $\pi = 3$) into a precise circular circumference (stage 2, $\pi_{susa} = 3.125$). Its reciprocal, \textit{57;36}, allowed inverse calculations from arcs back to linear hexagon-parts.
	
	The historical evolution of $\Xi$ is continuous and well documented:
	\begin{itemize}
		\item \textbf{Gudea (ca. 2100 BCE):} $\Xi_1 = 375/360 = \textit{1;2;30}$ (effective $\pi = 3.125$)
		\item \textbf{Ptolemy (ca. 150 CE):} $\Xi_2 = 377/360 = \textit{1;2;50}$ (effective $\pi = 3.141\bar{6}$)
		\item \textbf{Modern:} $\Xi = 376.9911\ldots/360$ (effective $\pi = 3.14159\ldots$)
	\end{itemize}
	The accuracy increased by only 0.5\% over four millennia, confirming that the conceptual framework was established early and underwent only incremental numerical refinement.
	
	\subsection*{$\Xi$ as a Babylonian Predecessor}
	
	The preceding sections have shown that this heading is clearly not an anachronism for several reasons already explained. Furthermore, this work does not suggest that the Babylonians invented a separate coefficient for calculating arc length and called it the radian. Due to the unique characteristics of the sexagesimal system, the single coefficient discussed here and called $\Xi$, can be applied to both functions simultaneously without any modification. In detail, the ancient Babylonian system used one-``sixth" of a circle (60 parts), and according to this innovative very effective thinking and system this also functioned precisely for one part of circle as of its complete 360 parts. Our modern system against that, uses three of these same``sixths" of circle (180 parts or degrees), which we then reduce back to one part and if necessary, multiplied by 60 to one- ``sixth" of a circle (hexagon). In essence, however, this is completely identical to the original ancient system. Everyone can decide which one is in detail the original, or the more efficient. 
	
	That is why the comparison with the modern radian measure leads to the surprising, but undeniable, result. The Babylonian coefficient $\Xi$ functions analogously to the modern radian measure: It is the proportionality factor between the number of parts (hexagon segments, functionally equivalent to modern degrees) and the arc length, represented by the stretching factor for each of these parts on the circumference.\\ 
	\\
	Specifically:
	This comparison serves only to demonstrate the similarity or technical equivalence of the coefficients. Their refinement over the last 4000 years can be determined using the Old Babylonian and modern values.
	\begin{itemize}
		\item $\Xi = \textit{1;2;30}$ (sexagesimal) Conceptually, this exactly corresponds to our modern $\pi/3 \approx 1.047197\ldots$ (modern) — the radian for $60^\circ$.
		\item $\Xi = \textit{0;1;2;30}$ (sexagesimal) in decimal $\Xi/60$ Conceptually, this exactly corresponds to our modern $\pi/180 \approx 0.017453\ldots$ — the radian for $1^\circ$.\\ While in modern this value applies to all decimal parts of 1 degree fraction, the Old Babylonian $\Xi$ also was applied without changes in sexagesimal notation (but decimal $\Xi /3600$ and so on), to the smaller fractions of 1 part, which we today call: arc minutes and arc seconds. 
		\item $1/\Xi = \textit{57;36}$ Conceptually, this exactly corresponds to our modern $180/\pi \approx 57.2958\ldots$ — the number of degrees at one (1) radian.
	\end{itemize}
	
	\subsection*{Implications for the History of Mathematics}
	
	Modern mathematics often views $\pi$ and radian measure as abstract constants of a modern era. This analysis reveals a seamless algorithmic evolution starting in Mesopotamia 4,000 years ago, leading through Greek antiquity (Ptolemy's chord tables) to modern analysis. The two-stage refinement process — first the hexagon, then the stretched circle — is not a later interpretation, but is built into the very structure of our modern $\pi$, which still contains the Old Babylonian components (Here simultaneously in the order of their introduction.): {\large $ 3 \cdot \Xi = \pi$.}
	
	The question of whether the Babylonians knew $\pi$ or $\Xi$ first is thus answered: $\Xi$ was the initially created coefficient to combine the hexagon with the circle — not as a mathematical abstraction, but as a metrological necessity for precise measurements on the sky and on the ground. $\pi$ was the logical consequence and result when high precision was needed, obtained by multiplying the old coefficient ``3'' with the stretching factor $\Xi$.\\
	\\
	This work certainly does not present all possible and historically used approaches on mathematical calculations involving circles around the world, but it does illustrate the Old Babylonian method, which can be traced through classical Greece up to our modern time, into our formulas and coefficients. We should be proud, beside that we have continuously refined these initial insights and achievements born out of practical necessity, we still continue to utilize their legacy in our daily life.
	\\
	The sexagesimal place-value system, with its inherent agnosticism toward the absolute magnitude of a number (60 = 1 = 1/60), allowed the same coefficient to serve simultaneously for the whole circle (360), for its sixth part (60), and for its 360th part (1). And exactly this is not an anachronistic equation of superficially similar numbers but a recognition of functional and historical continuity. The Babylonians lacked the abstract concept of an angle or measuring in radians, but they possessed its functional and conceptual equivalents: a proportional coefficient that related the number of parts (degrees) to the arc length. The refinements only demonstrate the development based on their clear defined predecessor and initial concept.
	\subsection*{Final Remarks}
	
	While, beside of this knowing, this analysis cannot uncover every nuance of the ancient Babylonian mind, and much work remains for future scholars, it is time to move beyond the notion that Babylonian mathematics was limited to superficial or merely practical. The sophisticated use of modular coefficients, the two-stage refinement process and the functional equivalence to the modern radian measure, all point to a level of geometric and astronomical sophistication that has been underestimated. It is not their knowledge, but perhaps our own understanding of their legacy, that remains incomplete.
	
	The goddess Nisaba, who in Gudea's dream held the golden stylus and the lapis-lazuli tablet of the starry heaven, demanded precise measures. The Babylonians delivered — and their coefficient $\Xi$ lives on in every calculation of arc length, every conversion from degrees to radians, and every use of $\pi$ today.\\
	\\
	\\
	Jens Kleb, Germany


\newpage
	\begin{thebibliography}{30}
		\bibitem{Beckmann}
		Beckmann, P. (1970). \textit{A History of $\pi$}.
		
		\bibitem{Bewer1919}
		Bewer, J. A. (1919). “The Cylinders of Gudea.” \textit{The American Journal of Semitic Languages and Literatures}, Vol. 35, No. 4, pp. 239–256.
		
		\bibitem{BIB-Orient}
		Bibliotheca Orientalis 21 (1964). ``Babylonische Mathematik," Nederlands Instituut voor het Nabije Oosten, pp. 45-52 by J.E. Hofmann.
		
		\bibitem{BruinsRutten1961}
		Bruins, E. M., \& Rutten, M. (1961). \textit{Textes mathématiques de Suse}. Paris: Paul Geuthner.
		
		\bibitem{Bruins1950}
		Bruins, E.M. (1950). Aperçu sur les mathématiques babyloniennes. \textit{Revue d'histoire des sciences et de leurs applications}, tome 3, no 4, pp. 301-314.
		
		\bibitem{CDLI_Gudea}
		CDLI (2025). No. P232300 / P232301 (Gudea Cylinders A and B). \url{https://cdli.earth/}.
		
		\bibitem{CDLI3}
		CDLI (2025). No. P254835. \url{https://cdli.earth/P254835} (TMS3 = MDP 34, 03 = SB13089).
		
		\bibitem{CDLI2}
		CDLI (2025). No. P254827. \url{https://cdli.earth/P254827} (TMS2 = MDP 34, 02 = SB13088).
		
		\bibitem{CDLI1}
		CDLI (2025). No. P254816. \url{https://cdli.earth/P254816} (TMS1 = MDP 34, 01 = SB13090).
		
		\bibitem{CDLI4}
		CDLI (2023). No. P255045. \url{https://cdli.earth/P255045} (YBC7243 = MCT136).
		
		\bibitem{CDLI5}
		CDLI (2023). No. P255026. \url{https://cdli.earth/P255026} (YBC5022 = MCT132).
		
		\bibitem{Edzard1997}
		Edzard, D. O. (1997). \textit{Gudea and His Dynasty}. Royal Inscriptions of Mesopotamia, Early Periods 3/1. Toronto: University of Toronto Press.
		
		\bibitem{ETCSL_c2.1.7}
		ETCSL (Electronic Text Corpus of Sumerian Literature). “The building of Ninĝirsu’s temple (Gudea, cylinders A and B)” – c.2.1.7. Oxford University. \url{https://etcsl.orinst.ox.ac.uk/}.
		
		\bibitem{Friberg2005}
		Friberg, J. (2005). \textit{Unexpected Links between Egyptian and Babylonian Mathematics}. Singapore: World Scientific.
		
		\bibitem{Hoyrup2002}
		Höyrup, J. (2002). \textit{Lengths, Widths, Surfaces: A Portrait of Old Babylonian Algebra and Its Kin}. New York: Springer.
		
		\bibitem{Kleb2021}
		Kleb, J. (2021). On the evidence of.... \textit{arXiv:2112.13379 [math.HO]}.
		
		\bibitem{Kleb2022}
		Kleb, J. (2022). Greek geometry theorems on Old Babylonian Cuneiform Tablets. \textit{Cuneiform Digital Library Initiative CDLP 25.0}.
		
		\bibitem{Kleb2023a}
		Kleb, J. (2023a). What is $\pi$? On the emergence of the circle constant. \textit{Conference Proceedings 2023}. \url{https://doi.org/10.6084/m9.figshare.24042636}.
		
		\bibitem{Kleb2023b}
		Kleb, J. (2023b). Himmelsbeobachtungen und deren Messung. \textit{Historia i Świat}, Tome 12, 209-236.
		
		\bibitem{Kleb2023c}
		Kleb, J. (2023c). The Reconstruction of the Hipparchan Dioptra, \textit{SSRN Digital Library: \url{http://dx.doi.org/10.2139/ssrn.4636206}}.
		
		\bibitem{LiuJian2016}
		Liu Jian (2016). “Sumerian Temple Building Ceremonies in Gudea Cylinders A and B.” \textit{Journal of Ancient Civilizations}, Vol. 31, pp. 1–27.
		
		\bibitem{Muroi}
		Muroi, K. (2016). The oldest example of $\pi$ in Sumer. \textit{arXiv:1610.03380 [math.HO]}.
		
		\bibitem{Neugebauer1945}
		Neugebauer, O. (1945). \textit{Mathematical Cuneiform Texts}. American Oriental Series Volume 29. New Haven: American Oriental Society.
		
		\bibitem{Neugebauer1954}
		Neugebauer, O. (1954). The exact Sciences in Antiquity [Review]. \textit{Revue d'histoire des sciences}, 7(3), 286–288.
		
		\bibitem{Neugebauer1969}
		Neugebauer, O. (1969). \textit{The exact sciences in antiquity} (2nd ed.). New York: Dover Publications.
		
		\bibitem{Neugebauer1962}
		Neugebauer, O. (1962). \textit{The astronomical tables of Al-Khwarizmi}. Kobenhavn.
		
		\bibitem{PtolemyAlmagest}
		Ptolemy, C. (ca. 150 CE). \textit{Almagest}. Book I, chapters 10–11.
		
		\bibitem{ORACC}
		Robson, E. (2014). The Digital Corpus of Cuneiform Mathematical Texts. \url{http://oracc.ub.uni-muenchen.de/dccmt/}
		
		\bibitem{Rey2020}
		Rey, S. (2020). “Gudea’s Temple Building and the Role of the Goddess Nisaba.” In \textit{Proceedings of the 65th Rencontre Assyriologique Internationale}, Paris.
		
		\bibitem{Roemer2010}
		Römer, W. H. Ph. (2010). \textit{Die Zylinderinschriften von Gudea}. Alter Orient und Altes Testament 376. Münster: Ugarit-Verlag.
		
		\bibitem{Rottlaender1979}
		Rottländer, R. C. A. (1979). \textit{Vormetrische Längenmaße}. Vienna: Verlag der Österreichischen Akademie der Wissenschaften.
		
		\bibitem{Rottlaenderwebsite}
		Rottländer, R. C. A. (2010). \textit{www.vormetrische-laengeneinheiten.de}
		
		\bibitem{ScribaSchreiber2005}
		Scriba, C. J., \& Schreiber, P. (2005). \textit{5000 Jahre Geometrie} (2nd ed.). Berlin: Springer. [English translation: \textit{5000 Years of Geometry}, Basel: Birkhäuser, 2015].
		
		\bibitem{Verderame2020}
		Verderame, L. (2020). “I sogni di Gudea: oniromanzia e politica nella Mesopotamia del III millennio a.C.” \textit{Studi Epigrafici e Linguistici}, Vol. 35–36, pp. 123–138.
		
		\bibitem{WorldHistory_Nisaba}
		World History Encyclopedia (2021). “Nisaba” – \url{https://www.worldhistory.org/Nisaba/}.
		
	\end{thebibliography}
\end{document}